\documentclass[10pt]{article}
\usepackage{amsfonts}

\usepackage{amscd}

\newtheorem{theorem}{Theorem}[section]
\newtheorem{corollary}[theorem]{Corollary}

\newtheorem{proposition}[theorem]{Proposition}

\newtheorem{definition}{Definition}[section]
\newtheorem{remark}{Remark}[section]

\begin{document}

\noindent {\Large Symplectic or Contact Structures on }\newline
{\Large Lie Groups } 
\[
\]

\noindent \textbf{Yu.Khakimdjanov and M.Goze}\medskip \newline
\textit{Universit\'{e} de Haute Alsace}\newline
\textit{4, rue des Fr\`{e}res Lumi\`{e}re, 68093 Mulhouse C\'{e}dex, France.}
\newline
\textit{E-mails: Y.Hakimjanov@uha.fr, \medskip M.Goze@uha.fr} \newline
\textbf{A.Medina} \medskip \newline
\textit{Universit\'{e} de Montpellier II\newline
case 051, 34095 Montpellier C\'{e}dex 5, France} \newline
\textit{E-mail: medina@math.univ-montp2.fr} 
\[
\]

\noindent {\LARGE Introduction}

In the sequel $G$ stands for a Lie group (supposed to be connected as a
matter of simplicity) with Lie algebra $\frak{g}:=T_{\varepsilon }\left(
G\right) $ , where $\varepsilon $ is the unit of $G$. If $G$ is endowed with
a left invariant differential $1$-form $\alpha ^{+}$ such that 
\[
\alpha ^{+}\wedge \left( d\alpha ^{+}\right) ^{p}\neq 0 
\]

\noindent where $2p+1$ is the dimension of $G,$ we will say that the pair $%
\left( G,\,\alpha ^{+}\right) $ is a contact Lie group and that $\left( 
\frak{g},\,\alpha \right) $ is a contact Lie algebra; here $\alpha :=\alpha
_{\varepsilon }^{+}$ .

Following Lichnerowicz-Medina \cite{LM} a pair $\left( G,\,\Omega
^{+}\right) $ where $\Omega ^{+}$ is a left invariant symplectic form, is
termed a symplectic Lie group and the corresponding infinitesimal object $%
\left( \frak{g},\,\omega \right) ,$ where $\omega :=\Omega _{\varepsilon }^{+},$ is referred to as a symplectic Lie algebra.

In \cite{MR} (see also \cite{DM1}, \cite{DM2}) a method of construction of
symplectic Lie algebras, called ''Symplectic Double Extension'', is
described. According to the theorem 2.5 in \cite{MR} every nilpotent
symplectic Lie algebra is obtained from a sequence of ''Symplectic Double
Extension'' starting from the trivial abelian Lie algebra consisting on only
one element.

This result immediately implies that every nilpotent contact Lie algebra can
be obtained by two operations, namely: the ''Symplectic Double Extension''
and the contactization.

Corresponding to those operations are inverse operations, well-known by
geometers, the symplectic reduction and the symplectization.

Here are some few words about the main results and the organization of this
work.

The section 1 gives a geometric description of the Contact Lie groups. In
the theorem 1, they arise to be fibre bundles with connections the fibre
being one dimensional, over a reductive homogeneous space 
\[
H\stackrel{i}{\hookrightarrow }G\stackrel{\pi }{\rightarrow }M=G/H 
\]

\noindent provided with a symplectic form, satisfying $\pi ^{*}\left( \Omega
\right) =\widehat{\Omega }$ where $\widehat{\Omega }$ is the curvature form
of the connection $\alpha ^{+}$ (see \cite{BW}). Sometimes $G/H$ is a
symplectic Lie group.

The section 2 supplies a necessary and sufficient condition for a filiform
Lie group to possess a left invariant contact form (see the Theorem 4). Such
a contact form is unique up to a non zero scalar multiple and has a simple
expression in terms of an adapted basis (see the Theorem 5).

Here it is convenient to recall some known facts. According to a result from
Gromov, every Lie group of odd dimension admits a non necessary left
invariant contact form. A symplectic Lie group $\left( G,\,\Omega
^{+}\right) $ is endowed with a left invariant affine structure (see \cite
{Chu}) defined by the following formulas for $a,\,b,\,c$ in $\frak{g}$ 
\[
\omega \left( ab,c\right) =-\omega \left( b,\left[ a,c\right] \right) 
\]
\[
\bigtriangledown _{a^{+}}b^{+}:=\left( ab\right) ^{+} 
\]

\noindent where $a^{+}$ is the left invariant vector field on $G$ such that $%
a_{\varepsilon }^{+}=a.$ Such connection $\bigtriangledown $ is fundamental
in the description of the symplectic Lie groups and specially
K\"{a}hl\'{e}rian Lie groups (see \cite{DM1}, \cite{DM2}). Unlike the
symplectic case there exists contact Lie groups with no left invariant
affine structure. This is what happens for semi-simple contact Lie groups.
More surprising, there even exists nilpotent contact Lie groups that never
admit such an affine structure: a direct verification allows us to check
that the example of Benoist (of dimension 11) supplied in \cite{B} is among
them.

Our work ends by supplying all nilpotent symplectic Lie algebras of
dimension $\leq 6.$

In this paper the following standard convention will be used without
explicit mentioning: for a concrete basis $X_{1},\ldots ,X_{n}$ of a Lie
algebra only those brackets $[X_{i},X_{j}]$ which are nonzero and for which $%
i<j$ will be explicitly defined .

\section{Contact Lie Groups as principal bundles with connection}

The aim of this section is to prove the more or less known following results
(see \cite{BW}, \cite{G1}, \cite{G2}).

\begin{theorem}
\textit{Let }$\left( G,\alpha ^{+}\right) $\textit{\ be a connected contact
Lie group and }$H$\textit{\ the isotropy subgroup of }$\alpha :=\alpha
_{\varepsilon }^{+}$ \textit{, for the coadjoint action. Then}

$\left( a\right) $\textit{\ The Lie group }$H$\textit{\ is 1-dimensional and
the homogeneous space }$M:=G/H$\textit{\ is reductive in the sense of Nomizu.%
}

$\left( b\right) $\textit{\ The form }$\alpha ^{+}$\textit{\ is a
''connection form'' on the canonical principal bundle} 
\begin{equation}
H\stackrel{i}{\hookrightarrow }G\stackrel{\pi }{\rightarrow }M=G/H  \label{1}
\end{equation}

\textit{the curvature form\ }$\widetilde{\Omega }$\textit{\ of which
satisfies the condition }$\widetilde{\Omega }=d\alpha ^{+}.$\textit{\ }

$\left( c\right) $\textit{\ There exists a symplectic form }$\Omega $\textit{%
\ on }$M$\textit{\ such that }$\pi ^{*}\left( \Omega \right) =\widetilde{%
\Omega }.$

$\left( d\right) $\textit{\ The canonical action of }$G$\textit{\ on }$%
\left( G/H_{0},\Omega _{0}\right) $\textit{\ is Hamiltonian, where }$H_{0}$%
\textit{\ is the connected component of the unit in }$H$\textit{\ and }$%
\Omega _{0}=p^{*}\left( \Omega \right) $\textit{\ with }$p$\textit{\ being
the natural projection of }$G/H_{0}$\textit{\ onto }$G/H.$
\end{theorem}

\noindent {\it Proof.}  It is clear that $H:=\left\{ \sigma \in G:\mbox{Ad}^{*}\left(
\sigma \right) \left( \alpha \right) =\alpha \right\} $ is a closed (hence
embedded) subgroup of $G,$ the Lie algebra $L\left( H\right) =\left\{ x\in 
\frak{g}:\mbox{ad}^{*}\left( x\right) \left( \alpha \right) =0\right\} $ of
which coincides with the radical Rad$\left( d\alpha \right) $ of the
bilinear form $d\alpha .$ Set $\dim G=2p+1$ and let's prove that $\dim H=1.$
As $\alpha ^{+}$ is a contact form, Ker$\alpha \cap $Rad$\left( d\alpha
\right) =\left\{ 0\right\} $ so that one has 
\[
\begin{array}{r}
0\leq \dim \left( \mbox{Ker}\alpha +\mbox{Rad}\left( d\alpha \right) \right)
=\dim \mbox{Ker}\alpha +\dim \mbox{Rad}\left( d\alpha \right)  \\ 
=2p+\dim \mbox{Rad}\left( d\alpha \right) \leq 2p+1
\end{array}
\]

\noindent that is $0\leq \dim H\leq 1$ \cite{G1}$.$

If $\dim H=0,$ a fortiori $G$ is of odd dimension, as the manifold $G/H$ and 
$Orb\left( \alpha \right) ,$ orbit of $\alpha $ via Ad$_{G}^{*},$ are
diffeomorphic. This is absurd. Thus $\dim H=1.$

Let's prove that $M:=G/H$ is reductive.. Let $z\in L\left( H\right) $ such
that $\alpha \left( z\right) =1;$ one has $L\left( H\right) =\Bbb{R}z.$ Set $%
\frak{m}:=\left\{ x\in \frak{g}:\alpha \left( x\right) =0\right\} ,$ then we
get $\frak{g}=L\left( H\right) \oplus \frak{m}.$ Furthermore for $x\in \frak{%
m}$ and $\tau \in H$ we have 
\[
\begin{tabular}{l}
$\alpha \left( \mbox{Ad}\left( \tau ^{-1}\right) \left( x\right) \right) =$Ad%
$^{*}\left( \tau \right) \left( \alpha \right) \left( x\right) =\alpha
\left( x\right) =0$%
\end{tabular}
\]

\noindent i.e. Ad$^{*}\left( H\right) \left( \frak{m}\right) \subset \frak{m}%
.$

Let $z^{+}$ be the left invariant vector field in $G$ with $z_{\varepsilon
}^{+}=z.$ For every $X\in T_{\sigma }\left( G\right) ,$ let 
\[
\theta _{\sigma }\left( X\right) :=\alpha _{\sigma }^{+}\left( X\right)
\,z_{\sigma }^{+} 
\]

\noindent Let's check that $\theta $ is a connection form.

Denote $x^{*}$ the vertical (relative to the fibration (\ref{1})) vector
field on $G$ associated to $x\in L\left( H\right) .$ For $\sigma \in G,$ one
has 
\[
x_{\sigma }^{*}:=\frac{d}{dt}\mid _{t=0}\left( \sigma \exp tx\right)
=x_{\sigma }^{+} 
\]

\noindent As $x=\lambda z$ for some $\lambda \in \Bbb{R},$ it follows 
\[
x_{\sigma }^{+}:=\left( L_{\sigma }\right) _{*,\varepsilon }\left( x\right)
=\left( L_{\sigma }\right) _{*,\varepsilon }\left( \lambda z\right) =\lambda
\left( L_{\sigma }\right) _{x,\varepsilon }\left( z\right) =\lambda
z_{\sigma }^{+}=x_{\sigma }^{*} 
\]

\noindent Thus 
\[
\theta _{\sigma }\left( x_{\sigma }^{*}\right) =\theta _{\sigma }\left(
\lambda z_{\sigma }^{+}\right) =\lambda \theta _{\sigma }\left( z_{\sigma
}^{+}\right) =\lambda \alpha _{\sigma }^{+}\left( z_{\sigma }^{+}\right)
z_{\sigma }^{+}=\lambda z_{\sigma }^{+}=x_{\sigma }^{*} 
\]

\noindent Now let's prove that for every $\tau \in H$ we have $\left(
R_{\tau }\right) ^{*}\theta =$Ad$\left( \tau ^{-1}\right) \theta .$ From the
equalities Ad$^{*}\left( \tau \right) \alpha =\alpha $ and $L_{\tau
}^{*}\alpha ^{+}=\alpha ^{+},$ its follows that $R_{\tau }^{*}\alpha
^{+}=\alpha ^{+}$ for every $\tau \in H.$

\noindent For $X_{\sigma }\in T_{\sigma }\left( G\right) $ and for every $%
\tau \in H$, $\sigma \in G$ we have$:$%
\[
\left( R_{\tau }^{*}\theta \right) \left( X_{\sigma }\right) =\theta
_{\sigma \tau }\left( \left( R_{\tau }\right) _{*,\alpha }\,X_{\sigma
}\right) =\alpha _{\sigma \tau }^{+}\left( \left( R_{\tau }\right)
_{*,\alpha }\,X_{\sigma }\right) z_{\sigma r}^{+} 
\]
\[
\theta _{\sigma }\left( X_{\sigma }\right) \,=\alpha _{\sigma }^{+}\left(
X_{\sigma }\right) z_{\sigma }^{+} 
\]

\noindent As $\alpha _{\sigma \tau }^{+}\left( \left( R_{\tau }\right)
_{*,\alpha }\,X_{\sigma }\right) =\alpha _{\sigma }^{+}\left( X_{\sigma
}\right) $ it follows that $R_{\tau }^{*}\theta =\theta .$ Hence we must
check that $\theta =$Ad$\left( \tau ^{-1}\right) \cdot \theta =R_{\tau
}^{*}\theta $ for $\tau \in H,$ that is 
\[
\theta \left( \left( R_{\tau }\right) _{*}X\right) =\mbox{Ad}\left( \tau
^{-1}\right) \theta \left( X\right) 
\]

\noindent for every $\tau \in H.$ But this arises from the fact that $H$ is
commutative. Thus $\theta $ is a connection form.

Let $\widetilde{\Omega }$ be the curvative form of $\theta .$ From the fact
that $\dim H=1,$ the relation 
\[
d\theta \left( X,Y\right) =-\frac{1}{2}\left[ \theta \left( X\right) ,\theta
\left( Y\right) \right] +\widetilde{\Omega }\left( X,Y\right) 
\]

\noindent for all $X,Y$ $\in $ $T_{\sigma }\left( G\right) $ then reads 
\[
d\theta \left( X,Y\right) =\widetilde{\Omega }\left( X,Y\right) =d\alpha
^{+}\left( X,Y\right) 
\]

\noindent Let's prove $\left( c\right) .$ As $\widetilde{\Omega }=d\alpha
^{+}$ we will have $L_{\sigma }^{*}\,\widetilde{\Omega }=\widetilde{\Omega }$
for all $\sigma \in G.$ Furthermore, for every $\tau $ in $H,$ 
\[
R_{\tau }^{*}\,\widetilde{\Omega }=R_{\tau }^{*}\left( d\alpha ^{+}\right)
=d\left( R_{\tau }^{*}\alpha ^{+}\right) =d\alpha ^{+}=\widetilde{\Omega } 
\]

\noindent Let $\left[ \sigma \right] \in M$ and $u,\,v$ in $T_{\left[ \sigma
\right] }\left( M\right) .$ Set 
\[
\Omega _{\left[ \sigma \right] }\left( u,v\right) :=\widetilde{\Omega }%
_{\sigma }\left( u_{\sigma },v_{\sigma }\right) 
\]

\noindent where $u_{\sigma }$ (respectively $v_{\sigma }$) is the horizontal
lifts of $u$ (respectively of $v$) at $\sigma .$ Let's see first that $%
\Omega $ is well defined. Let $u_{\sigma \tau },$ $v_{\sigma \tau }$ be the
horizontal lifts of $u$ and $v$ at $\sigma \tau $ with $\tau \in H.$ One has 
\[
\widetilde{\Omega }_{\sigma \tau }\left( u_{\sigma \tau },v_{\sigma \tau
}\right) =\widetilde{\Omega }_{\sigma \tau }\left( \left( R_{\tau }\right)
_{*,\alpha }\,u_{\sigma },\left( R_{\tau }\right) _{*,\alpha }\,v_{\sigma
}\right) =\widetilde{\Omega }_{\sigma }\left( u_{\sigma },v_{\sigma }\right) 
\]

\noindent as $\widetilde{\Omega }$ is $R_{\tau }$ invariant, $\tau \in H.$

In addition the equalities 
\[
\pi ^{*}\left( d\Omega \right) =d\left( \pi ^{*}\Omega \right) =d\widetilde{%
\Omega }=d\left( d\alpha ^{+}\right) =0 
\]

\noindent imply that $d\Omega =0$ and taking into account the following 
\[
\pi ^{*}\left( \Omega ^{p}\right) =\left( \pi ^{*}\Omega \right) ^{p}=(%
\widetilde{\Omega })^{p}=\left( d\alpha ^{+}\right) ^{p}\neq 0 
\]

\noindent we then deduce that $\Omega $ is symplectic and invariant by the
canonical action of $G$ on $M.$

The canonical map $p:G/H_{0}\rightarrow G/H,\quad \sigma H_{0}\rightarrow
\sigma H,$ is obviously a covering map. Let $\Omega _{0}:=p^{*}\left( \Omega
\right) .$ It is clear that $\Omega _{0}$ is symplectic and invariant by the
canonical action of $G$ on $G/H_{0}.$ Furthermore, one has $R_{\tau }^{*}\,%
\widetilde{\Omega }=\widetilde{\Omega }$ for all $\tau $ in $H_{0}.$ We have 
$\pi _{0}^{*}\Omega _{0}=$ $\widetilde{\Omega }$ where $\pi
_{0}:G\rightarrow G/H_{0}=:M_{0}$ is the canonical injection.

We are going to prove now that the canonical action 
\[
\phi :G\times M_{0}\rightarrow M_{0}\qquad \left( \sigma ,\left[ \rho
\right] \right) \mapsto \left[ \sigma \rho \right] =:\phi _{\sigma }\left(
\left[ \rho \right] \right) 
\]

\noindent is a Hamiltonian action. The action $\phi $ is symplectic. For $%
x\in \frak{g}$ , let $\left[ \widetilde{x}\right] $ be the fundamental
vector field on $M_{0}$ associated to $x.$ It is clear that the following
diagram is commutative: 
\[
\begin{tabular}{lllll}
& $\quad \;G$ & $\;\stackrel{L_{\sigma }}{\longrightarrow }\;$ & $G$ &  \\ 
& $\pi _{0}\downarrow $ &  & $\downarrow \pi _{0}$ &  \\ 
& $\quad \;M_{0}$ & $\stackrel{\phi _{\sigma }}{\;\longrightarrow }\;$ & $%
M_{0}$ & 
\end{tabular}
\]

\noindent for every $\sigma $ in $G.$ Let's denote by $\widetilde{x}$ the
horizontal lift (relative to $\theta _{0}$) of $\left[ \widetilde{x}\right] $
on the total space of the fiber with connection $H_{0}\hookrightarrow
G\rightarrow M_{0}.$ This vector field is invariant under the $R_{\tau }$
for $\tau \in H_{0}.$ Moreover, as the flows of $\widetilde{x}$ and $\left[ 
\widetilde{x}\right] $ are the same via $\pi _{0}$ and 
\[
\left[ \widetilde{x}\right] _{\left[ \sigma \right] }:=\frac{d}{dt}\left|
_{t=0}\right. \exp tx\cdot \left[ \sigma \right] =\frac{d}{dt}\left|
_{t=0}\right. [\exp tx\cdot \sigma ] 
\]

\noindent it follows that the flow of $\widetilde{x}$ consists on left
translations on $G.$ Thus $\widetilde{x}$ is a right invariant vector field
on $G.$ Let's emphasize on the fact that we are not pretending that $%
\widetilde{x}$ is the right invariant vector field $x^{-}$ associated to $x.$
We only have $\widetilde{x}=y^{-}$ for some $y\in \frak{g}$ satisfying $\pi
_{*,\varepsilon }\left( y\right) =\left[ x\right] _{\pi _{0}\left(
\varepsilon \right) }.$

Let $\left( \varphi _{t}\right) _{t}$ be the flow of $\widetilde{x}=y^{-}.$
One has 
\[
0=\mathcal{L}\left( y^{-}\right) \alpha ^{+}=\left( d\circ i\left(
y^{-}\right) +i\left( y^{-}\right) \circ d\right) \alpha ^{+}=d\left( \alpha
^{+}(y^{-})\right) +(i\left( y^{-}\right) \circ d)\left( \alpha ^{+}\right) 
\]

\noindent where $\mathcal{L}$ is the derivative. But one also has 
\[
i\left( y^{-}\right) d\alpha ^{+}=i\left( y^{-}\right) \pi _{0}^{*}\left(
\Omega _{0}\right) =\pi _{0}^{*}\left( i\left[ \widetilde{x}\right] (\Omega
_{0})\right) 
\]

\noindent so that 
\begin{equation}
0=d\left( \alpha ^{+}\left( y^{-}\right) \right) +\pi _{0}^{*}\left( i\left[ 
\widetilde{x}\right] \Omega _{0}\right)  \label{2}
\end{equation}

\noindent Let's consider the function $f_{y}:G\rightarrow \Bbb{R}$, \quad $%
f_{y}\left( \sigma \right) :=\alpha _{\sigma }^{+}\left( y_{\sigma
}^{-}\right) .$ Let's prove that $f_{y}$ can be projected by $\pi _{0}.$ Let 
$Y\in T_{\sigma }\left( G\right) .$ From $($\ref{2}$)$ we have, 
\[
d\left( \alpha ^{+}\left( y^{-}\right) \right) \left( Y_{\sigma }\right)
=Y_{\sigma }\left( \alpha ^{+}\left( y^{-}\right) \right) =Y_{\sigma }\left(
f_{y}\right) =-\Omega _{0}\left( \left[ \widetilde{x}\right] ,\left( \pi
_{0}\right) _{*,\sigma }Y_{\sigma }\right) 
\]

\noindent Consequently, if $Y$ is tangent to the fiber (that is if $\left(
\pi _{0}\right) _{*,\sigma }Y_{\sigma }=0$), we'll have $Y_{\sigma }\left(
f_{y}\right) =0$ for every $\sigma $ in $G.$ Hence $f$ is constant along $%
H_{0}.$ This implies the existence of a smooth function $J_{y}:M_{0}%
\rightarrow \Bbb{R}$ such that $J_{y}\circ \pi _{0}=f_{y}.$%

\bigskip

The following result is a complement of the theorem. It is directly proved
by taking into account the ideas provided in the proof of the theorem.

\begin{corollary}
\cite{G1}

$\left( a\right) $ If $\left( G,\alpha ^{+}\right) $ is a contact Lie group
of non discrete center $Z\left( G\right) $, then the quotient Lie group $%
G/Z\left( G\right) $ has left invariant symplectic form $\Omega ^{+}$ such
that $\pi ^{*}\Omega ^{+}=-d\alpha ^{+},$ where $\pi :G\rightarrow G/Z\left(
G\right) $ is the canonical projection.

$\left( b\right) $ Conversely, if $\left( K,\Omega ^{+}\right) $ is a
symplectic Lie group, every Lie group $G$ with Lie algebra $L\left( G\right)
:=\Bbb{R}_{\omega }{\times }L\left( K\right) $ (the central
extension of $L\left( K\right) $ by $\Bbb{R}$ via $\omega $), where $\omega
:=\Omega _{\varepsilon }^{+},$ admits a left invariant contact form $\alpha
^{+}$ satisfying $\pi ^{*}\Omega ^{+}=-d\alpha ^{+}$ (that is $\pi
^{*}\omega =-d\alpha $).
\end{corollary}

\begin{remark}
$\left( a\right) $ Notice that the manifold $G/H$ in the above theorem can
be identified with the orbit $Orb\left( \alpha \right) $ of $\alpha \in 
\frak{g}^{*}$ for the coadjoint representation.

$\left( b\right) $ If two contact Lie algebras $\left( \frak{g}_{1},\alpha
_{1}\right) $ and $\left( \frak{g}_{2},\alpha _{2}\right) $ with non trivial
centers are contacto-isomorphic (that is there exists an isomorphism of Lie
algebras $\varphi :\frak{g}_{1}\rightarrow \frak{g}_{2},$ such that, $%
\varphi ^{*}\left( \alpha _{2}\right) =\alpha _{1}$), then the symplectic
Lie algebras $\left( \frak{g}_{1}/Z\left( \frak{g}_{1}\right) ,\omega
_{1}\right) $ et $\left( \frak{g}_{2}/Z\left( \frak{g}_{2}\right) ,\omega
_{2}\right) $ are symplecto-isomorphic (that is there exists an isomorphism
of Lie algebras $\varphi :\frak{g}_{1}/Z\left( \frak{g}_{1}\right)
\rightarrow \frak{g}_{2}/Z\left( \frak{g}_{2}\right) ,$ such that, $\varphi
^{*}\left( \omega _{2}\right) =\omega _{1}$) and $\pi _{i}^{*}\omega
_{i}=-d\alpha _{i},\quad i=1,2;$ where $\pi _{i}:\frak{g}_{i}\rightarrow 
\frak{g}_{i}/Z\left( \frak{g}_{i}\right) $ are the canonical projections $.$

$\left( c\right) $ Let $\eta ^{+}$ be the left invariant $\frak{g}$-valued
invariant1-form on $G$ defined by $\eta _{\varepsilon }\left( x\right) =x$
for every $x\in \frak{g}.$ As $M:=G/H$ and $M_{0}:=G/H_{0}$ are reductive as
stated in theorem, then the $L\left( H\right) $ is a component of $\eta $
gives rise to a connection on the fiber bundles $H\stackrel{i}{%
\hookrightarrow }G\stackrel{\pi }{\rightarrow }M$ and $H\stackrel{i}{%
_{0}\hookrightarrow }G\stackrel{\pi _{0}}{\rightarrow }M_{0}$ which are
invariant under the action of $G$ (\cite{KN}, page 103). Such connections
coincide with the ones described in the above theorem.
\end{remark}

\section{Contact or Symplectic Filiform Lie algebras}

\subsection{Filiform Lie algebras (basic definitions and results)}

Let $\frak{g}$ be a nilpotent Lie algebra of dimension $n$.\ Let

\[
C^{0}\frak{g}\supset C^{1}\frak{g}\supset ...\supset C^{n-2}\frak{g}\supset
C^{n-1}\frak{g}=\{0\} 
\]

\noindent be the central descending series of $\frak{g}$, where $C^{0}\frak{g%
}=\frak{g}$, $C^{i}\frak{g}=\left[ \frak{g},C^{i-1}\frak{g}\right] $, $1\leq
i\leq n-1$.

\begin{definition}
A Lie algebra $\frak{g}$ of dimension $\geq 3$ is called \textit{filiform}
if $\dim C^{k}\frak{g}=n-k-1$ for $k=1,...,n-1$.
\end{definition}

We remark that the filiform Lie algebras have the maximal possible nilindex,
that is $n-1$. These algebras are the ''least'' nilpotent.\medskip

\noindent \textbf{Examples of filiform Lie algebras}

For each $n\in \Bbb{N}$ there exists several $\left( n+1\right) $%
-dimensional filiform Lie algebras which are specially remarkable. In the
following description, the brackets are given relative to a basis $\left(
X_{0},X_{1},...,X_{n}\right) .$

\begin{enumerate}
\item  The Lie algebra $L_{n}:$

It is the simplest $\left( n+1\right) $-dimensional filiform Lie algebra.
Its non trivial brackets are given by: 
\[
\left[ X_{0},X_{i}\right] =X_{i+1},\quad i=1,...,n-1.
\]

\item  The Lie algebra $Q_{n}\quad (n=2k+1):$%
\[
\begin{array}{l}
\left[ X_{0},X_{i}\right] =X_{i+1},\quad i=1,...,n-1 \\ 
\left[ X_{i},X_{n-i}\right] =\left( -1\right) ^{i}X_{n},\quad i=1,...,k
\end{array}
\]

In the basis $\left( Z_{0},Z_{1},...,Z_{n}\right) ,$ where $%
Z_{0}=X_{0}+X_{1},\quad Z_{i}=X_{i},\quad i=1,...,n;$ this Lie algebra is
defined by 
\[
\begin{array}{l}
\left[ Z_{0},Z_{i}\right] =Z_{i+1},\quad i=1,...,n-2 \\ 
\left[ Z_{i},Z_{n-i}\right] =\left( -1\right) ^{i}Z_{n},\quad i=1,...,k
\end{array}
\]

\item  The Lie algebra $R_{n}:$%
\[
\begin{array}{l}
\left[ X_{0},X_{i}\right] =X_{i+1},\quad i=1,...,n-1 \\ 
\left[ X_{1},X_{j}\right] =X_{j+2},\quad j=2,...,n-2
\end{array}
\]

\item  The Lie algebra $W_{n}:$%
\[
\begin{array}{l}
\left[ X_{0},X_{i}\right] =X_{i+1},\quad i=1,...,n-1 \\ 
\left[ X_{i},X_{j}\right] =\frac{6\left( i-1\right) !\left( j-1\right)
!\left( j-i\right) }{\left( i+j\right) !}X_{i+j+1},\quad 1\leq i,j\leq
n-2,\quad i+j+1\leq n
\end{array}
\]
This Lie algebra can be defined also relative to a basis $\left(
Y_{1},Y_{2},...,Y_{n+1}\right) $ by the brackets 
\[
\left[ Y_{i},Y_{j}\right] =\left( j-i\right) Y_{i+j},\quad i+j\leq n+1
\]

\item  \medskip The Lie algebra $T_{n}\quad \left( n=2k\right) :$%
\[
\begin{array}{l}
\left[ X_{0},X_{i}\right] =X_{i+1},\quad i=1,...,n-1 \\ 
\left[ X_{k-i-1},X_{k+i}\right] =\left( -1\right) ^{i}X_{n},\quad
i=0,1,...,k-2
\end{array}
\]

\item  The Lie algebra $T_{n}\quad \left( n=2k+1\right) :$%
\[
\begin{array}{l}
\left[ X_{0},X_{i}\right] =X_{i+1},\quad i=1,...,n-1 \\ 
\left[ X_{k-i-1},X_{k+i+j}\right] =\left( -1\right)
^{i}C_{i+j}^{i}X_{n+j-1},\quad i=0,1,...,k-2,\quad j=0,1
\end{array}
\]

\item  The Lie algebra $P_{n}\quad \left( n=2k\right) $%
\[
\begin{array}{l}
\left[ X_{0},X_{i}\right] =X_{i+1},\quad i=1,...,n-1;\quad \left[
X_{k-1},X_{k}\right] =X_{n}; \\ 
\left[ X_{k-i-1},X_{k+i}\right] =\left( -1\right) ^{i}\left( 1-\frac{2}{%
\left( k-1\right) \left( k-2\right) }C_{i+1}^{i-1}\right) X_{n},\quad
i=1,...,k-2; \\ 
\left[ X_{k-i-2},X_{k+i+j-1}\right] =\left( -1\right) ^{i}\frac{2}{\left(
k-1\right) \left( k-2\right) }C_{i+j}^{i}X_{n+j-2},\;0\leq i\leq k-3;\;j=0,1.
\end{array}
\]
\end{enumerate}

Let $\frak{g}$ be a $m$-dimensional filiform Lie algebra. It is naturally
filtered by descending central series and we can associate to $\frak{g}$ a
graded Lie algebra $gr\,\frak{g}$ which is also filiform. This Lie algebra
is defined on the vector space 
\[
gr\,\frak{g}=\oplus _{i=1}^{m-1}\frak{g}_{i} 
\]

\noindent where $\frak{g}_{i}=C^{i-1}\frak{g}/C^{i}\frak{g}$, by the
brackets $\left[ x+C^{i}\frak{g},y+C^{j}\frak{g}\right] =\left[ x,y\right]
+C^{i+j}\frak{g},$ $x\in C^{i-1}\frak{g}$, $y\in C^{j-1}\frak{g}$.

\begin{proposition}
\cite{Ve} Let $\frak{g}$ be a $m$-dimensional filiform Lie algebra. Then the
graded Lie algebra $gr\,\frak{g}$ is isomorphic to $L_{m-1}$ , if $m$ is
odd, and isomorphic to $L_{m-1}$ or $Q_{m-1}$, if $m$ is even.
\end{proposition}

Let $\Delta $ be the set of pairs of integers $\left( k,r\right) $ such that 
$1\leq k\leq n-1$, $2k+1<r\leq n$, $r\geq 4$ (if $n$ is odd we suppose that $%
\Delta $ contain also the pair $\left( \frac{n-1}{2},n\right) $).\ For any
element $(k,r)\in \Delta $, we can associate the 2-cocycle for the Chevalley
cohomology of $L_{n}$ with coefficients in the adjoint module denoted $\Psi
_{k,r}$ and defined by

\[
\Psi _{k,r}\left( X_{i},X_{j}\right) =-\Psi _{k,r}\left( X_{j},X_{i}\right)
=\left( -1\right) ^{k-i}C_{j-k-1}^{k-i}X_{i+j+r-2k-1} 
\]

\noindent \noindent if $1\leq i\leq k<j\leq n$ , $i+j+r-2k-1\leq n$ and $%
\Psi _{k,r}\left( X_{i},X_{j}\right) =0$ otherwise.\ We remark that this
formula for $\Psi _{k,r}$ is uniquely determined from the conditions :

\[
\Psi _{k,r}\left( X_{k},X_{k+1}\right) =X_{r} 
\]

\[
\Psi _{k,r}\left( X_{i},X_{j}\right) \in Z^{2}\left( L_{n},L_{n}\right) 
\]

\begin{proposition}
\cite{Ve}Any $(n+1)$-dimensional filiform Lie algebra law $\mu \in F_{m}$ is
isomorphic to $\mu _{0}+\Psi $ where $\mu _{0}$ is the law of $L_{n}$ and $%
\Psi $ is a 2-cocycle defined by 
\[
\Psi ={\sum_{{\left( k,r\right) \in \Delta }}}a_{k,r}\Psi _{k,r}
\]
and verifying the relation $\Psi \circ \Psi =0$ with

\[
\Psi \circ \Psi \left( x,y,z\right) =\Psi \left( \Psi \left( x,y\right)
,z\right) +\Psi \left( \Psi \left( y,z\right) ,x\right) +\Psi \left( \Psi
\left( z,x\right) ,y\right) 
\]
\end{proposition}

\begin{definition}
Let $\frak{g}$ be a $(n+1)$-dimensional filiform Lie algebra with law $\mu $%
. A basis $\left( X_{0},X_{1},...,X_{n}\right) $ of $\frak{g}$ is called 
\textit{adapted}, if $\left[ X_{i},X_{j}\right] =\mu _{0}\left(
X_{i},X_{j}\right) +\Psi \left( X_{i},X_{j}\right) ,\quad 0\leq i,j\leq n.$
\end{definition}

\begin{proposition}
Let $\frak{g}$ be a filiform Lie algebra of dimension $\geq 4.$ Then $Der%
\frak{g}$ is solvable.
\end{proposition}

\noindent {\it Proof}. 
Consider an adapted basis $\left( X_{0},X_{1},...,X_{n}\right) $ of $\frak{g}%
.$ As the central descending series of $\frak{g}$ is an invariant flag under
all derivations it is sufficient to show that the ideal $\left\langle
X_{1},...,X_{n}\right\rangle $ is also an invariant. Let $d\in Der\frak{g}$
and $d\left( X_{1}\right) =\sum_{i=0}^{n}a_{i}X_{i}.$ For the $4$%
-dimensional filiform Lie algebra $\overline{\frak{g}}=\frak{g}/C^{3}\frak{g}
$ (it is isomorphic to $L_{3}$) we have the derivation $\overline{d}$ with $%
\overline{d}\left( \overline{X}_{1}\right) =\sum_{i=0}^{3}a_{i}\overline{X}%
_{i}.$ This is possible only if $a_{0}=0.$

\bigskip

Let $\frak{g}$ be a Lie algebra. Consider in $Der\frak{g}$ a maximally
abelian subalgebra $\frak{t}$ consisting of semisimple endomorphisms (a such
subalgebra is called torus of $\frak{g}$). According to a theorem by Mostow 
\cite{Mo} two such subalgebras are conjugated by an inner automorphism. The
common dimension of a tori on $\frak{g}$ is called \textit{rank }of $\frak{g}
$. Note that for nilpotent $\frak{g}$ the rank cannot exceed the codimension
of the derived ideal since $\frak{g}$ is generated by any vector subspace of 
$\frak{g}$ complementary to the derived ideal. For a filiform Lie algebra
the only possible ranks are 0, 1 and 2.

\begin{proposition}
\cite{GK2}Let $\frak{g}$ be a filiform Lie algebra of dimension $n+1$ and of
rank 2. Then $\frak{g}$ is isomorphic to $L_{n}$ if $n$ is even and
isomorphic to $L_{n}$ or $Q_{n}$ if $n$ is odd.
\end{proposition}

The following theorem gives a description of the filiform Lie algebras of
rank 1.

\begin{theorem}
\cite{GK2}Let $\frak{g}$ be a filiform Lie algebra of dimension $n+1\geq 7$
and of rank 1.\ There is a basis $\left( Y_{0},Y_{1},...,Y_{n}\right) $ of $%
\frak{g}$ such that $\frak{g}$ is one of the following families of Lie
algebras: \medskip 

$\left( i\right) 
\begin{tabular}[b]{l}
$
\begin{array}{l}
\frak{g}=A_{n+1}^{r}\left( \alpha _{1},...,\alpha _{t}\right) ,\quad 1\leq
r\leq n-3,\quad t=\left[ \frac{n-r-1}{2}\right] , \\ 
\left[ Y_{0},Y_{i}\right] =Y_{i+1},\quad 1\leq i\leq n-1, \\ 
\left[ Y_{i},Y_{j}\right] =\left( \sum_{k=i}^{t}\alpha _{k}\left( -1\right)
^{k-i}C_{j-k-1}^{k-i}\right) Y_{i+j+r},\;1\leq i\leq j\leq n,\;i+j+r\leq n.
\end{array}
$%
\end{tabular}
$\medskip 

$\left( ii\right) 
\begin{tabular}{l}
$
\begin{array}{l}
\frak{g}=B_{n+1}^{r}\left( \alpha _{1},...,\alpha _{t}\right)
,\;n=2m+1,\;1\leq r\leq n-4,\;t=\left[ \frac{n-r-2}{2}\right] , \\ 
\left[ Y_{0},Y_{i}\right] =Y_{i+1},\quad 1\leq i\leq n-2, \\ 
\left[ Y_{i},Y_{n-i}\right] =(-1)^{i}Y_{n},\quad 1\leq i\leq m, \\ 
\left[ Y_{i},Y_{j}\right] =\left( \sum_{k=i}^{t}\alpha _{k}\left( -1\right)
^{k-i}C_{j-k-1}^{k-i}\right) Y_{i+j+r}, \\ 
1\leq i<j\leq n-1,\quad i+j+r\leq n-1.
\end{array}
$%
\end{tabular}
$\medskip 

$\left( iii\right) 
\begin{tabular}{l}
$
\begin{array}{l}
\frak{g}=C_{n+1}\left( \alpha _{1},...,\alpha _{t}\right) ,\quad
n=2m+1,\quad t=m-1, \\ 
\left[ Y_{0},Y_{i}\right] =Y_{i+1},\quad 1\leq i\leq n-2, \\ 
\left[ Y_{i},Y_{n-i}\right] =(-1)^{i}Y_{n},\quad 1\leq i\leq m, \\ 
\left[ Y_{i},Y_{n-i-2k}\right] =(-1)^{i}\alpha _{k}Y_{n},\quad 1\leq k\leq
m-1,\quad 1\leq i\leq n-2k-1
\end{array}
$%
\end{tabular}
$\medskip 

\noindent where $C_{q}^{s}$ are the binomial coefficients (we suppose that $%
C_{q}^{s}=0$ if $q<0$ or $q<s$), $\left( \alpha _{1},...,\alpha _{t}\right) $
are the parameters satisfying the polynomial relations emanating from
Jacobi's identity and at least one parameter $\alpha _{i}\neq 0$. A maximal
torus of derivations is spanned by $d$, where :

If $\frak{g}=A_{n+1}^{r}\left( \alpha _{1},...,\alpha _{t}\right) $\textbf{\
: } 
\[
d\left( Y_{0}\right) =Y_{0},\quad d\left( Y_{i}\right) =(i+r)Y_{i},\quad
1\leq i\leq n.
\]

If $\frak{g}=B_{n+1}^{r}\left( \alpha _{1},...,\alpha _{t}\right) $\textbf{: 
} 
\[
d\left( Y_{0}\right) =Y_{0},\quad d\left( Y_{i}\right) =(i+r)Y_{i},\quad
1\leq i\leq n-1,\quad d\left( Y_{n}\right) =\left( n+2r\right) Y_{n}.
\]

If \textbf{: }$\frak{g}=C_{n+1}\left( \alpha _{1},...,\alpha _{t}\right) $: 
\[
d\left( Y_{0}\right) =0,\quad d\left( Y_{i}\right) =Y_{i},\quad 1\leq i\leq
n-1,\quad d\left( Y_{n}\right) =2Y_{n}.
\]
\end{theorem}

\begin{remark}
$\left( a\right) $ Let $n\geq 13$ and $r=1.$ Then, up to isomorphism, there
are only four Lie algebras $\frak{g}$ of rank 1 if $n$ is even and three Lie
algebras of rank 1 if $n$ is odd \cite{KH}, \cite{KH1}: If $n$ is even and $%
\geq 14,$ then $\frak{g}$ is isomorphic to one of the Lie algebras $%
R_{n},\,W_{n},\,T_{n},\,P_{n}.$ If $n$ is odd and $\geq 13$, $\frak{g}$ is
isomorphic to one of the Lie algebras $\,R_{n},\,W_{n},\,T_{n}.$ If $n=12,$
then $\frak{g}$ is isomorphic to one of the Lie algebras $R_{12},$ $W_{12},$ 
$T_{12}.$

$\left( b\right) $ The laws $C_{n+1}\left( \alpha _{1},...,\alpha
_{t}\right) $ satisfy the Jacobi's identity for all values of parameters $%
\left( \alpha _{1},...,\alpha _{t}\right) $.

$\left( c\right) $ Let $\frak{g}$ be a Lie algebra belonging to one of the
families $(i),(ii),(iii)$ and at least one of parameters $\alpha _{i}$ be
different to zero.\ Then we can transform one of these parameters to 1 using
the automorphism $\psi $ defined by $\psi \left( X_{0}\right) =aX_{0}$, $%
\psi \left( X_{1}\right) =bX_{1}$ (this is a unique type of automorphisms
preserving the torus and the property of basis to be adapted). Modulo this
transformation we have a classification up to isomorphism of filiform Lie
algebras of rank 1.
\end{remark}

\subsection{\noindent Symplectization and contactization of the Filiform Lie
algebras}

The following result shows that the class of Filiform Lie algebras is closed
respect to the contactization and symplectization process described in the
section 1.

\begin{theorem}
\textit{Let} $(G,\alpha ^{+})$ \textit{be a contact filiform Lie group. Then
the quotient }$G/Z(G)$\textit{\ is a symplectic filiform Lie group.
Conversely, if }$(K,\omega )$\textit{\ is a symplectic filiform Lie group,
then every central extension} 
\[
0\rightarrow \Bbb{R\rightarrow }G\rightarrow K\rightarrow 0
\]
\textit{following }$\omega $\textit{\ is a contact filiform Lie group.}
\end{theorem}

\noindent  {\it Proof.}
Let $\frak{g}=L\left( G\right) $ be the Lie algebra of $G.$ For an adapted
basis $(X_{0},X_{1},\ldots ,X_{2p})$ of $\frak{g}$ we have $\left[
X_{0},X_{i}\right] =X_{i+1},\quad i=1,\ldots ,2p-1$ and $Z(\frak{g})=\Bbb{R}%
\cdot X_{2p}.$ Following the corollary of theorem 1 the quotient $\frak{g}/Z(%
\frak{g})$ is a symplectic Lie algebra. Let $\pi :\frak{g\rightarrow }%
\widetilde{\frak{g}}=\frak{g}/Z(\frak{g})$ be the canonical projection. Then
we have 
\[
\left[ \pi (X_{0}),\pi (X_{i})\right] =\pi (X_{i+1}),\quad i=1,\ldots ,2p-2; 
\]

\noindent and the Lie algebra $\widetilde{\frak{g}}$ is also filiform.

Conversely suppose that $(\widetilde{\frak{g}},\omega )$ is a symplectic
filiform Lie algebra of dimension $2p.$ Consider an adapted basis $%
(Y_{0},Y_{1},\ldots ,Y_{2p-1})$ of $\widetilde{\frak{g}}.$ As $\omega $ is a
non degenerated form, there exists $0\leq k\leq 2p-2,$ such that $\omega
(Y_{k},Y_{2p-1})\neq 0.$ Let $\frak{g}=\widetilde{\frak{g}}\oplus _{\omega }%
\Bbb{R}$ be the central extension following $\omega .$ The central
descending sequence $\left\{ C^{i}\frak{g}\right\} $ satisfies the condition 
$\dim C^{i-1}\frak{g}/C^{i}\frak{g}=1$ for all $2\leq i\leq 2p-2$ because
this property holds in $\widetilde{\frak{g}}.$ As $\omega
(Y_{k},Y_{2p-1})\neq 0$ we have also $\dim C^{2p-2}\frak{g}/C^{2p-1}\frak{g}%
=1$ and $\dim C^{2p-1}\frak{g}=1.$ Thus the nilindex of $\frak{g}$ is equal
to $2p$ and $\frak{g}$ is filiform.

\subsection{Existence of a left invariant contact form}

\begin{theorem}
\textit{Let }$G$\textit{\ be a }$\left( 2p+1\right) $\textit{-dimensional
filiform Lie group and }$\frak{g}$\textit{\ its Lie algebra. Suppose that
the law }$\mu $\textit{\ of }$\frak{g}$\textit{\ is written in an adapted
basis by the formula} 
\[
\mu =\mu _{0}+\sum_{(k,r)\in \Delta }a_{k,r}\Psi _{k,r}
\]
\textit{Then }$G$\textit{\ admits a left invariant contact form if and only
if} 
\[
A_{j}:=\sum_{s=0}^{j-1}(-1)^{s}a_{p-j+s,2p-2(j-s-1)}C_{2j-s-2}^{s}\neq
0,\quad j=1,2,\ldots ,p-1;
\]

\textit{if this property holds the linear form }$\alpha =b_{0}\alpha
_{0}+b_{1}\alpha _{1}+\ldots +b_{2p}\alpha _{2p}$\textit{\ is a contact form
on }$\frak{g}$\textit{\ if and only if }$b_{2p}\neq 0.$
\end{theorem}

\noindent {\it Proof. }
\textbf{\ }Let $A_{j}\neq 0,\quad j=1,2,\ldots ,p-1.$ We have

\[
\begin{array}{ll}
\left[ X_0,X_{2p-1}\right] & =X_{2p} \\ 
\left[ X_1,X_{2p-2}\right] & =\sum_{k=1}^{p-1}a_{k,2k+2}\Psi
_{k,2k+2}(X_1,X_{2p-2})=A_{p-1}X_{2p} \\ 
\left[ X_2,X_{2p-3}\right] & =\sum_{k=2}^{p-1}a_{k,2k+2}\Psi
_{k,2k+2}(X_2,X_{2p-3})=A_{p-2}X_{2p} \\ 
& \ldots \\ 
\left[ X_{p-1},X_p\right] & =a_{p-1,2p}\Psi _{p-1,2p}(X_{p-1},X_p)=A_1X_{2p}
\end{array}
\]

\noindent and $\left[ X_{j},X_{m}\right] =0,$ if $m\geq 2p-j.$ In the dual
basis $\left( \alpha _{0},\alpha _{1},\ldots ,\alpha _{2p}\right) $ of $%
\frak{g}^{*}$ the previous brackets give 
\begin{eqnarray*}
d\left( \alpha \right) &=&-\alpha _{0}\wedge \alpha _{2p-1}-A_{p-1}\alpha
_{1}\wedge \alpha _{2p-2}-A_{p-2}\alpha _{2}\wedge \alpha _{2p-3}-\ldots \\
&&-A_{1}\alpha _{p-1}\wedge \alpha _{p}+\sum_{m<2p-j-1}b_{j,m}\alpha
_{j}\wedge \alpha _{m}
\end{eqnarray*}

\noindent where $\alpha =a_{0}\alpha _{0}+a_{1}\alpha _{1}+\ldots
+a_{2p-1}\alpha _{2p-1}+\alpha _{2p}.$

\noindent Thus 
\[
\left( d\alpha \right) ^{p}=\left( -1\right) ^{p}p!A_{1}A_{2}\cdots
A_{p-1}\alpha _{0}\wedge \alpha _{1}\wedge \alpha _{2}\wedge \cdots \wedge
\alpha _{2p-1} 
\]

\noindent and $\alpha \wedge \left( d\alpha \right) ^{p}\neq 0.$ This means
that $\alpha $ is a contact form.

Conversely, we suppose now that the Lie algebra $\frak{g}$ admits a contact
form $\alpha .$ We put 
\[
\alpha =b_{0}\alpha _{0}+b_{1}\alpha _{1}+\ldots +b_{2p}\alpha _{2p}. 
\]

\noindent As $Z\left( \frak{g}\right) $ is not included in Ker$\alpha ,$
then $b_{2p}\neq 0.$ We have 
\[
d\alpha = 
\begin{array}[t]{l}
b_{0}(d\alpha _{0})+b_{1}(d\alpha _{1})+\ldots +b_{2p}(d\alpha _{2p})= \\ 
=b_{2}(-\alpha _{0}\wedge \alpha _{1})+b_{3}(-\alpha _{0}\wedge \alpha
_{2})+b_{4}(-\alpha _{0}\wedge \alpha _{3}-a_{1,4}\alpha _{1}\wedge \alpha
_{2})+\ldots + \\ 
+b_{j}(\sum_{l+s<j}t_{l,s}\alpha _{l}\wedge \alpha _{s})+\ldots + \\ 
+b_{2p}(-\alpha _{0}\wedge \alpha _{2p-1}-A_{p-1}\alpha _{1}\wedge \alpha
_{2p-2}-A_{p-2}\alpha _{2}\wedge \alpha _{2p-3}-\ldots - \\ 
-A_{1}\alpha _{p-1}\wedge \alpha _{p}+\sum_{m<2p-j-1}b_{j,m}\alpha
_{j}\wedge \alpha _{m}).
\end{array}
\qquad 
\begin{array}[t]{l}
\\ 
\\ 
\left( \ast \right)
\end{array}
\]
$\qquad $

\noindent As the basis $\left( X_{0},X_{1},...,X_{2p}\right) $ is adapted we
have $X_{2p}\perp d\alpha _{i},\quad 0\leq i\leq 2p,$ and thus $X_{2p}\perp
d\alpha .$ But $\left( d\alpha _{2p}\right) ^{p}\neq 0$ and 
\[
(d\alpha )^{p}=\lambda \alpha _{0}\wedge \alpha _{1}\wedge \alpha _{2}\wedge
\cdots \wedge \alpha _{2p-1} 
\]

\noindent with $\lambda \neq 0.$ Let us examine the terms of the expression $%
\left( *\right) $ which appear in the non null product $\left( d\alpha
_{2p}\right) ^{p}$. In this expression, only there is one term containing
the form $\alpha _{2p-1};$ it is the term $-b_{2p}\alpha _{0}\wedge \alpha
_{2p-1}.$ We deduce that 
\[
(d\alpha )^{p}=-b_{2p}\alpha _{0}\wedge \alpha _{2p-1}\wedge \theta 
\]
where $\theta \in \wedge ^{2p-2}\frak{g}$, $\theta \neq 0$ and $\alpha
_{0},\alpha _{2p-1}$ does not appear in $\theta .$ Let as examine now $%
\alpha _{2p-2}.$ As this containing in $(d\alpha )^{p},$ then $\theta =$ $%
-b_{2p}A_{p-1}\alpha _{1}\wedge \alpha _{2p-2}\wedge \theta _{1}$ with $%
\theta _{1}\in \wedge ^{2p-4}\frak{g}$, $\theta _{1}\neq 0.$ In fact in the
expression $t_{l,s}\alpha _{l}\wedge \alpha _{s}$ with $l+s<j,\quad j<2p,$
the index $s$ cannot be equal to $2p-2$ except the case $l=0.$ Likewise in
the term $b_{j,m}\alpha _{j}\wedge \alpha _{m}$ , we have $m\neq 2p-2$ if $%
j\neq 0.$

Let us suppose now that the terms 
\[
-b_{2p}\alpha _{0}\wedge \alpha _{2p-1},\ -b_{2p}A_{p-1}\alpha _{1}\wedge
\alpha _{2p-2},\ \ldots \ ,\ -b_{2p}A_{p-k}\alpha _{k}\wedge \alpha
_{2p-k-1} 
\]

\noindent of the expression $\left( *\right) $ are the factors of $(d\alpha
)^{p}.$ Then $A_{p-1},\ldots ,A_{p-k}\neq 0$ for $1\leq k<p-1.$ In the same
way we show that the term $-b_{2p}A_{p-k}\alpha _{k+1}\wedge \alpha
_{2p-k-2} $ is also a factor of $(d\alpha )^{p}.$ By induction we have 
\[
(d\alpha )^{p}=(-1)^{p}p!b_{2p}^{p}A_{1}A_{2}\cdots A_{p-1}\alpha _{0}\wedge
\alpha _{1}\wedge \alpha _{2}\wedge \cdots \wedge \alpha _{2p-1}\neq 0 
\]

\noindent and $A_{1},A_{2},\ldots ,A_{p-1}\neq 0.$

\subsection{Classes of contacto-isomorphisms}

Two contact Lie algebras $\left( \frak{g}_{1},\alpha _{1}\right) $ and $%
\left( \frak{g}_{2},\alpha _{2}\right) $ called contacto-isomorphic if there
exists an isomorphism of Lie algebras $\varphi :\frak{g}_{1}\rightarrow 
\frak{g}_{2},$ such that, $\varphi ^{*}\left( \alpha _{2}\right) =\alpha
_{1}.$ The following result gives the classification up
contacto-isomorphisms of the contact forms on a filiform Lie algebras.

\begin{theorem}
\textit{Let }$\frak{g}$\textit{\ be a filiform }$\left( 2p+1\right) $\textit{%
-dimensional Lie algebra. Let us consider an adapted basis }$%
(X_{0},X_{1},\ldots ,X_{2p})$\textit{\ of }$\frak{g}$\textit{\ and its dual
basis }$(\alpha _{0},\alpha _{1},\ldots ,\alpha _{2p}).$\textit{\ If }$%
\alpha =a_{0}\alpha _{0}+a_{1}\alpha _{1}+\ldots +a_{2p}\alpha _{2p}$\textit{%
\ is a contact form on }$\frak{g}$\textit{\ , then the form }$\beta
=a_{2p}\alpha _{2p}$\textit{\ is also a contact form on }$\frak{g}$\textit{\
and }$(\frak{g}$\textit{,} $\alpha )$ \textit{is contacto-isomorphic to} $(%
\frak{g}$, $\beta ).$
\end{theorem}

\noindent {\it Proof. }
The fact that $\beta $ is a contact form on $\frak{g}$ is a consequence of
the proof of theorem 4. To prove the theorem it is sufficient to find an
automorphism $\varphi \in $Aut $\frak{g}$ such that $\varphi ^{*}\alpha
=\beta .$

Consider the derivation $d=-a_{2p-1}$ad $X_{0}.$ As ad $X_{0}(X_{i})=X_{i+1},
$ $i=1,\ldots ,2p-1,$ the automorphism $A=\exp d$ satisfies the following
property: 
\[
A(X_{2p-1})=X_{2p-1}-a_{2p-1}X_{2p},\quad A(X_{2p})=X_{2p}.
\]

\noindent Then 
\[
A^{*}\alpha _{2p}=a_{2p}\alpha _{2p}-a_{2p-1}\alpha
_{2p-1}+\sum_{i<2p-1}c_{i}\alpha _{i} 
\]

\noindent and 
\[
A^{*}\alpha =a_{2p}\alpha _{2p}+\sum_{j\leq 2p-2}q_{j}\alpha _{j} 
\]

\noindent We suppose now that 
\[
\alpha =a_{2p}\alpha _{2p}+\sum_{j\leq 2p-k}b_{j}\alpha _{j},\quad 2\leq
k\leq 2p-1, 
\]

\noindent and we prove the existence of an automorphism $\varphi \in $Aut $\frak{g}$ such that 
\[
\varphi ^{*}\alpha =a_{2p}\alpha _{2p}+\sum_{j\leq 2p-k-1}b_{j}^{\prime
}\alpha _{j}. 
\]

\noindent Consider the derivation $d=-\lambda b_{2p-k}$ad $X_{k-1}.$ From
the proof of the theorem 4 we have 
\[
\mbox{ad}X_{k-1}(X_{2p-k})=A_{p-k+1}X_{2p},\quad 2\leq k\leq p, 
\]

\noindent and 
\[
\mbox{ad}X_{k-1}(X_{2p-k})=-A_{k-p}X_{2p},\mbox{ si }k>p. 
\]

\noindent We have also ad$X_{k-1}(X_{i})=0$ for all index $i>2p-k.$ Let us
put 
\[
\lambda =\left\{ 
\begin{tabular}{lll}
$\frac{1}{A_{p-k+1}},$ & si & $2\leq k\leq p$ \\ 
$\frac{-1}{A_{k-p}},$ & si & $i>2p-2$%
\end{tabular}
\right. 
\]

\noindent Then the automorphism $\varphi =\exp d$ satisfies the required
condition. By induction we deduce the theorem.

\bigskip

\begin{remark}
$\left( a\right) $\textbf{\ }The theorem 5 only concerns the filiform case.
But, from a direct verification, we can affirm that this result remains true
for every nilpotent Lie algebra of dimension less or equal to 7.

$\left( b\right) $\textbf{\ }Suppose that $\alpha _{2p}$ is a contact form
on the filiform Lie algebra $\frak{g}$ . In general the contact Lie algebras 
$(\frak{g},\,\,a\alpha _{2p})$ and $(\frak{g},\,b\alpha _{2p})$ with $a\neq b
$ are not contacto-isomorphic. But these contact Lie algebras are always
contacto-isomorphic if $\frak{g}$ admits a semisimple derivation (see \cite
{GK2} and the subsection 2.1 for their description).

$\left( c\right) $\textbf{\ }If $K=\Bbb{C},$ then the theorems 3, 4, 5 and
the point $\left( a\right) $ of this remark are valid. But the point $\left(
b\right) $ of this remark must be modified (see \cite{GJK2}).
\end{remark}

\section{Symplectic Lie algebras of dimension $\leq 6$}

The studied relation between the classes of symplectic Lie algebras of
dimension $2p$ and the contact Lie algebras of dimension $2p+1$ and the
description of contact structures on a filiform Lie algebras and on a
nilpotent Lie algebras of dimension $\leq 7$ permits to obtain some
classification results about symplectic Lie algebras. The following theorem
gives a complete classification up to simplecto-isomorphism of the
symplectic Lie algebras of dimension $\leq 6.$

\begin{theorem}
\textit{Every nilpotent symplectic Lie algebra of the dimension }$\leq 6$%
\textit{\ is symplecto-isomorphic to one and only one of the following
symplectic Lie algebras. \medskip }
\end{theorem}

\textbf{Dimension 2}

\begin{enumerate}
\item  $\qquad 
\begin{array}[t]{l}
\Bbb{R}^{2},\medskip  \\ 
\omega =\alpha _{1}\wedge \alpha _{2}
\end{array}
$\medskip 
\end{enumerate}

\textbf{Dimension 4}

\begin{enumerate}
\item  \qquad $
\begin{array}[t]{l}
\left[ X_{1},X_{2}\right] =X_{3},\qquad \left[ X_{1},X_{3}\right]
=X_{4},\medskip  \\ 
\omega =\alpha _{1}\wedge \alpha _{4}+\alpha _{2}\wedge \alpha _{3}
\end{array}
\qquad $

\item  \qquad $
\begin{array}[t]{l}
\left[ X_{1},X_{2}\right] =X_{3},\medskip  \\ 
\omega =\alpha _{1}\wedge \alpha _{3}+\alpha _{2}\wedge \alpha _{4}
\end{array}
\qquad $

\item  \qquad $
\begin{array}[t]{l}
\Bbb{R}^{4},\medskip  \\ 
\omega =\alpha _{1}\wedge \alpha _{4}+\alpha _{2}\wedge \alpha _{3}
\end{array}
$\medskip 
\end{enumerate}

\textbf{Dimension 6}

\begin{enumerate}
\item  $\qquad 
\begin{array}[t]{l}
\left[ X_{1},X_{2}\right] =X_{3},\qquad \left[ X_{1},X_{3}\right]
=X_{4},\qquad \left[ X_{1},X_{4}\right] =X_{5}, \\ 
\left[ X_{1},X_{5}\right] =X_{6},\qquad \left[ X_{2},X_{3}\right]
=X_{5},\qquad \left[ X_{2},X_{4}\right] =X_{6},\medskip  \\ 
\omega =\alpha _{1}\wedge \alpha _{6}+(1-\lambda )\alpha _{2}\wedge \alpha
_{5}+\lambda \alpha _{3}\wedge \alpha _{4},\qquad \lambda \in \Bbb{R}%
\setminus \{0,1\}
\end{array}
$

\item  $\qquad 
\begin{array}[t]{l}
\left[ X_{1},X_{2}\right] =X_{3},\qquad \left[ X_{1},X_{3}\right]
=X_{4},\qquad \left[ X_{1},X_{4}\right] =X_{5}, \\ 
\left[ X_{1},X_{5}\right] =X_{6},\qquad \left[ X_{2},X_{3}\right]
=X_{6},\medskip \qquad  \\ 
\omega (\lambda )=\lambda (\alpha _{1}\wedge \alpha _{6}+\alpha _{2}\wedge
\alpha _{4}+\alpha _{3}\wedge \alpha _{4}-\alpha _{2}\wedge \alpha
_{5}),\qquad \lambda \in \Bbb{R}\setminus \{0\}
\end{array}
$

\item  $\qquad 
\begin{array}[t]{l}
\left[ X_{1},X_{2}\right] =X_{3},\qquad \left[ X_{1},X_{3}\right]
=X_{4},\qquad \left[ X_{1},X_{4}\right] =X_{5}, \\ 
\left[ X_{1},X_{5}\right] =X_{6},\medskip \qquad  \\ 
\omega =\alpha _{1}\wedge \alpha _{6}-\alpha _{2}\wedge \alpha _{5}+\alpha
_{3}\wedge \alpha _{4}\qquad 
\end{array}
$

\item  $\qquad 
\begin{array}[t]{l}
\left[ X_{1},X_{2}\right] =X_{3},\qquad \left[ X_{1},X_{3}\right]
=X_{4},\qquad \left[ X_{1},X_{4}\right] =X_{6}, \\ 
\left[ X_{2},X_{3}\right] =X_{5},\qquad \left[ X_{2},X_{5}\right]
=X_{6},\medskip  \\ 
\omega (\lambda _{1},\lambda _{2})=\lambda _{1}\alpha _{1}\wedge \alpha
_{4}+\lambda _{2}(\alpha _{1}\wedge \alpha _{5}+\alpha _{1}\wedge \alpha
_{6}+\alpha _{2}\wedge \alpha _{4}+\alpha _{3}\wedge \alpha _{5}),\qquad  \\ 
\lambda _{1}\in \Bbb{R},\quad \lambda _{2}\in \Bbb{R}\setminus \{0\}
\end{array}
$

\item  $\qquad 
\begin{array}[t]{l}
\left[ X_{1},X_{2}\right] =X_{3},\qquad \left[ X_{1},X_{3}\right]
=X_{4},\qquad \left[ X_{1},X_{4}\right] =-X_{6}, \\ 
\left[ X_{2},X_{3}\right] =X_{5},\qquad \left[ X_{2},X_{5}\right]
=X_{6},\medskip  \\ 
\begin{array}[t]{r}
\omega _{1}(\lambda _{1},\lambda _{2})=\lambda _{1}\alpha _{1}\wedge \alpha
_{4}+\lambda _{2}(\alpha _{1}\wedge \alpha _{5}+\alpha _{1}\wedge \alpha
_{6}+\alpha _{2}\wedge \alpha _{4}+\alpha _{3}\wedge \alpha _{5}), \\ 
\lambda _{1}\in \Bbb{R},\quad \lambda _{2}\in \Bbb{R}\setminus \{0\}
\end{array}
\qquad  \\ 
\begin{array}[t]{r}
\omega _{2}(\lambda )=\lambda (-\alpha _{1}\wedge \alpha _{6}+\alpha
_{3}\wedge \alpha _{4}+\frac{1}{2}\alpha _{1}\wedge \alpha _{4}+\frac{1}{4}%
\alpha _{1}\wedge \alpha _{5}+\frac{1}{4}\alpha _{2}\wedge \alpha _{4} \\ 
-\alpha _{3}\wedge \alpha _{5}),\quad \lambda \in \Bbb{R}\setminus \{0\}
\end{array}
\\ 
\begin{array}[t]{r}
\omega _{3}(\lambda )=\lambda (-\alpha _{1}\wedge \alpha _{6}+\alpha
_{3}\wedge \alpha _{4}+\frac{3}{2}\alpha _{1}\wedge \alpha _{4}-\frac{3}{4}%
\alpha _{1}\wedge \alpha _{5}+\frac{1}{4}\alpha _{2}\wedge \alpha _{4} \\ 
-\alpha _{3}\wedge \alpha _{5}),\quad \lambda \in \Bbb{R}\setminus \{0\}
\end{array}
\quad  \\ 
\begin{array}[t]{r}
\omega _{4}(\lambda )=\lambda (-\alpha _{1}\wedge \alpha _{6}+\alpha
_{3}\wedge \alpha _{4}-\frac{1}{2}\alpha _{1}\wedge \alpha _{4}+\frac{5}{4}%
\alpha _{1}\wedge \alpha _{5}+\frac{1}{4}\alpha _{2}\wedge \alpha _{4} \\ 
-\alpha _{3}\wedge \alpha _{5}),\quad \lambda \in \Bbb{R}\setminus \{0\}
\end{array}
\end{array}
$

\item  $\qquad 
\begin{array}[t]{l}
\left[ X_{1},X_{2}\right] =X_{3},\qquad \left[ X_{1},X_{3}\right]
=X_{4},\qquad \left[ X_{1},X_{4}\right] =X_{5}, \\ 
\left[ X_{2},X_{3}\right] =X_{6},\medskip \qquad  \\ 
\omega _{1}=\alpha _{1}\wedge \alpha _{6}+\alpha _{2}\wedge \alpha
_{4}+\alpha _{2}\wedge \alpha _{5}-\alpha _{3}\wedge \alpha _{4}, \\ 
\omega _{2}=-\alpha _{1}\wedge \alpha _{6}-\alpha _{2}\wedge \alpha
_{4}-\alpha _{2}\wedge \alpha _{5}+\alpha _{3}\wedge \alpha _{4}
\end{array}
$

\item  $\qquad 
\begin{array}[t]{l}
\left[ X_{1},X_{2}\right] =X_{4},\qquad \left[ X_{1},X_{4}\right]
=X_{5},\qquad \left[ X_{1},X_{5}\right] =X_{6}, \\ 
\left[ X_{2},X_{3}\right] =X_{6},\qquad \left[ X_{2},X_{4}\right]
=X_{6},\medskip  \\ 
\omega _{1}(\lambda )=\lambda (\alpha _{1}\wedge \alpha _{3}+\alpha
_{2}\wedge \alpha _{6}+\alpha _{4}\wedge \alpha _{5}),\quad \lambda \in \Bbb{%
R}\setminus \{0\} \\ 
\omega _{2}(\lambda )=\lambda (\alpha _{1}\wedge \alpha _{6}+\alpha
_{2}\wedge \alpha _{5}-\alpha _{3}\wedge \alpha _{4}),\quad \lambda \in \Bbb{%
R}\setminus \{0\}
\end{array}
$

\item  $\qquad 
\begin{array}[t]{l}
\left[ X_{1},X_{2}\right] =X_{4},\qquad \left[ X_{1},X_{4}\right]
=X_{5},\qquad \left[ X_{1},X_{5}\right] =X_{6}, \\ 
\left[ X_{2},X_{3}\right] =X_{6},\qquad \left[ X_{2},X_{4}\right]
=X_{6},\medskip  \\ 
\omega =\alpha _{1}\wedge \alpha _{6}+\alpha _{2}\wedge \alpha _{5}-\alpha
_{3}\wedge \alpha _{4}
\end{array}
\qquad $

\item  $\qquad 
\begin{array}[t]{l}
\left[ X_{1},X_{2}\right] =X_{4},\qquad \left[ X_{1},X_{4}\right]
=X_{5},\qquad \left[ X_{1},X_{5}\right] =X_{6}, \\ 
\left[ X_{2},X_{3}\right] =X_{6},\medskip \qquad  \\ 
\omega (\lambda )=\lambda (\alpha _{1}\wedge \alpha _{3}+\alpha _{2}\wedge
\alpha _{6}+\alpha _{4}\wedge \alpha _{5}),\quad \lambda \in \Bbb{R}^{+}
\end{array}
$

\item  $\qquad 
\begin{array}[t]{l}
\left[ X_{1},X_{2}\right] =X_{4},\qquad \left[ X_{1},X_{4}\right]
=X_{5},\qquad \left[ X_{1},X_{3}\right] =X_{6}, \\ 
\left[ X_{2},X_{4}\right] =X_{6},\medskip \qquad  \\ 
\omega _{1}=\alpha _{1}\wedge \alpha _{6}+\alpha _{2}\wedge \alpha
_{5}-\alpha _{2}\wedge \alpha _{6}-\alpha _{3}\wedge \alpha _{4}, \\ 
\omega _{2}=-\alpha _{1}\wedge \alpha _{6}-\alpha _{2}\wedge \alpha
_{5}+\alpha _{2}\wedge \alpha _{6}+\alpha _{3}\wedge \alpha _{4}
\end{array}
$

\item  $\qquad 
\begin{array}[t]{l}
\left[ X_{1},X_{2}\right] =X_{4},\qquad \left[ X_{1},X_{4}\right]
=X_{5},\qquad  \\ 
\left[ X_{2},X_{3}\right] =X_{6},\qquad \left[ X_{2},X_{4}\right]
=X_{6},\medskip  \\ 
\omega _{1}(\lambda )=\alpha _{1}\wedge \alpha _{6}+\alpha _{2}\wedge \alpha
_{5}+\lambda \alpha _{2}\wedge \alpha _{6}-\alpha _{3}\wedge \alpha
_{4},\quad \lambda \in \Bbb{R} \\ 
\omega _{2}(\lambda )=-\alpha _{1}\wedge \alpha _{6}-\alpha _{2}\wedge
\alpha _{5}+\lambda \alpha _{2}\wedge \alpha _{6}+\alpha _{3}\wedge \alpha
_{4},\quad \lambda \in \Bbb{R}
\end{array}
$

\item  $\qquad 
\begin{array}[t]{l}
\left[ X_{1},X_{2}\right] =X_{4},\qquad \left[ X_{1},X_{4}\right]
=X_{5},\qquad \left[ X_{1},X_{3}\right] =X_{6}, \\ 
\left[ X_{2},X_{3}\right] =-X_{5},\qquad \left[ X_{2},X_{4}\right]
=X_{6},\medskip  \\ 
\begin{array}[t]{r}
\omega (\lambda _{1},\lambda _{2})=\lambda _{1}\alpha _{1}\wedge \alpha
_{5}+\lambda _{2}\alpha _{2}\wedge \alpha _{6}+\left( \lambda _{1}+1\right)
\alpha _{3}\wedge \alpha _{4}, \\ 
\lambda _{1}\in \Bbb{R}\setminus \{0,-1\},\quad \lambda \in \Bbb{R}^{+}
\end{array}
\qquad 
\end{array}
$

\item  $\qquad 
\begin{array}[t]{l}
\left[ X_{1},X_{2}\right] =X_{4},\qquad \left[ X_{1},X_{3}\right]
=X_{5},\qquad \left[ X_{1},X_{4}\right] =X_{6}, \\ 
\left[ X_{2},X_{3}\right] =X_{6},\medskip  \\ 
\omega _{1}(\lambda )=\alpha _{1}\wedge \alpha _{6}+\lambda \alpha
_{2}\wedge \alpha _{5}+\left( \lambda -1\right) \alpha _{3}\wedge \alpha
_{4},\quad \lambda \in \Bbb{R}\setminus \{0,1\} \\ 
\omega _{2}(\lambda )=\alpha _{1}\wedge \alpha _{6}+\lambda \alpha
_{2}\wedge \alpha _{4}+\alpha _{2}\wedge \alpha _{5}+\alpha _{3}\wedge
\alpha _{5},\quad \lambda \in \Bbb{R}\setminus \{0\} \\ 
\omega _{3}=\alpha _{1}\wedge \alpha _{6}+\alpha _{2}\wedge \alpha _{4}+%
\frac{1}{2}\alpha _{2}\wedge \alpha _{5}-\frac{1}{2}\alpha _{3}\wedge \alpha
_{4}
\end{array}
$

\item  $\qquad 
\begin{array}[t]{l}
\left[ X_{1},X_{2}\right] =X_{4},\qquad \left[ X_{1},X_{4}\right]
=X_{6},\qquad \left[ X_{1},X_{3}\right] =X_{5},\medskip  \\ 
\omega _{1}=\alpha _{1}\wedge \alpha _{6}+\alpha _{2}\wedge \alpha
_{4}+\alpha _{3}\wedge \alpha _{5}, \\ 
\omega _{2}=\alpha _{1}\wedge \alpha _{6}-\alpha _{2}\wedge \alpha
_{4}+\alpha _{3}\wedge \alpha _{5},\quad  \\ 
\omega _{3}=\alpha _{1}\wedge \alpha _{6}+\alpha _{2}\wedge \alpha
_{5}+\alpha _{3}\wedge \alpha _{4}
\end{array}
$

\item  $\qquad 
\begin{array}[t]{l}
\left[ X_{1},X_{2}\right] =X_{4},\qquad \left[ X_{1},X_{4}\right]
=X_{6},\qquad \left[ X_{2},X_{3}\right] =X_{5},\medskip  \\ 
\omega _{1}=-\alpha _{1}\wedge \alpha _{5}+\alpha _{1}\wedge \alpha
_{6}+\alpha _{2}\wedge \alpha _{5}+\alpha _{3}\wedge \alpha _{4}, \\ 
\omega _{2}=\alpha _{1}\wedge \alpha _{5}-\alpha _{1}\wedge \alpha
_{6}-\alpha _{2}\wedge \alpha _{5}-\alpha _{3}\wedge \alpha _{4},\quad  \\ 
\omega _{3}=\alpha _{1}\wedge \alpha _{6}+\alpha _{2}\wedge \alpha
_{4}+\alpha _{3}\wedge \alpha _{5}, \\ 
\omega _{4}=\alpha _{1}\wedge \alpha _{6}+\alpha _{2}\wedge \alpha
_{5}-\alpha _{3}\wedge \alpha _{4}, \\ 
\omega _{5}=-\alpha _{1}\wedge \alpha _{6}-\alpha _{2}\wedge \alpha
_{5}+\alpha _{3}\wedge \alpha _{4}
\end{array}
$

\item  $\qquad 
\begin{array}[t]{l}
\left[ X_{1},X_{2}\right] =X_{5},\qquad \left[ X_{1},X_{3}\right]
=X_{6},\qquad  \\ 
\left[ X_{2},X_{4}\right] =X_{6},\qquad \left[ X_{3},X_{4}\right]
=-X_{5},\medskip  \\ 
\omega _{1}=\alpha _{1}\wedge \alpha _{6}+\alpha _{2}\wedge \alpha
_{3}-\alpha _{4}\wedge \alpha _{5},\quad  \\ 
\omega _{2}=\alpha _{1}\wedge \alpha _{6}-\alpha _{2}\wedge \alpha
_{3}+\alpha _{4}\wedge \alpha _{5}
\end{array}
$

\item  $\qquad 
\begin{array}[t]{l}
\left[ X_{1},X_{3}\right] =X_{5},\qquad \left[ X_{1},X_{4}\right]
=X_{6},\qquad \left[ X_{2},X_{3}\right] =X_{6},\medskip  \\ 
\omega =\alpha _{1}\wedge \alpha _{6}+\alpha _{2}\wedge \alpha _{5}+\alpha
_{3}\wedge \alpha _{4}
\end{array}
$

\item  $\qquad 
\begin{array}[t]{l}
\left[ X_{1},X_{2}\right] =X_{4},\qquad \left[ X_{1},X_{3}\right]
=X_{5},\qquad \left[ X_{2},X_{3}\right] =X_{6},\medskip  \\ 
\omega _{1}(\lambda )=\alpha _{1}\wedge \alpha _{6}+\lambda \alpha
_{2}\wedge \alpha _{5}+\left( \lambda -1\right) \alpha _{3}\wedge \alpha
_{4},\quad \lambda \in \Bbb{R}\setminus \{0,1\}, \\ 
\omega _{2}(\lambda )=\alpha _{1}\wedge \alpha _{5}+\lambda \alpha
_{1}\wedge \alpha _{6}-\lambda \alpha _{2}\wedge \alpha _{5}+\alpha
_{2}\wedge \alpha _{6}-2\lambda \alpha _{3}\wedge \alpha _{4}, \\ 
\qquad \quad \qquad \qquad \qquad \qquad \qquad \qquad \qquad \qquad \qquad
\qquad \lambda \in \Bbb{R}\setminus \{0\}, \\ 
\omega _{3}=-2\alpha _{1}\wedge \alpha _{6}+\alpha _{2}\wedge \alpha
_{4}-\alpha _{2}\wedge \alpha _{5}+\alpha _{3}\wedge \alpha _{4}
\end{array}
$

\item  $\qquad 
\begin{array}[t]{l}
\left[ X_{1},X_{2}\right] =X_{4},\qquad \left[ X_{1},X_{4}\right]
=X_{5},\qquad \left[ X_{1},X_{5}\right] =X_{6},\medskip  \\ 
\omega =\alpha _{1}\wedge \alpha _{3}+\alpha _{2}\wedge \alpha _{6}+\alpha
_{4}\wedge \alpha _{5}
\end{array}
$

\item  $\qquad 
\begin{array}[t]{l}
\left[ X_{1},X_{2}\right] =X_{3},\qquad \left[ X_{1},X_{3}\right]
=X_{4},\qquad  \\ 
\left[ X_{1},X_{4}\right] =X_{5},\qquad \left[ X_{2},X_{3}\right]
=X_{5},\medskip  \\ 
\omega _{1}=\alpha _{1}\wedge \alpha _{6}+\alpha _{2}\wedge \alpha
_{5}-\alpha _{3}\wedge \alpha _{4},\quad  \\ 
\omega _{2}=-\alpha _{1}\wedge \alpha _{6}-\alpha _{2}\wedge \alpha
_{5}+\alpha _{3}\wedge \alpha _{4}
\end{array}
$

\item  $\qquad 
\begin{array}[t]{l}
\left[ X_{1},X_{2}\right] =X_{4},\qquad \left[ X_{1},X_{4}\right]
=X_{6},\qquad \left[ X_{2},X_{3}\right] =X_{6}, \\ 
\\ 
\omega =\alpha _{1}\wedge \alpha _{5}+\alpha _{2}\wedge \alpha _{4}-\alpha
_{3}\wedge \alpha _{4}-\alpha _{3}\wedge \alpha _{5}
\end{array}
$

\item  $\qquad 
\begin{array}[t]{l}
\left[ X_{1},X_{2}\right] =X_{5},\qquad \left[ X_{1},X_{5}\right] =X_{6}, \\ 
\\ 
\omega =\alpha _{1}\wedge \alpha _{6}+\alpha _{2}\wedge \alpha _{5}+\alpha
_{3}\wedge \alpha _{4}
\end{array}
$

\item  $\qquad 
\begin{array}[t]{l}
\left[ X_{1},X_{2}\right] =X_{5},\qquad \left[ X_{1},X_{3}\right] =X_{6}, \\ 
\\ 
\omega _{1}=\alpha _{1}\wedge \alpha _{6}+\alpha _{2}\wedge \alpha
_{5}+\alpha _{3}\wedge \alpha _{4} \\ 
\omega _{2}=\alpha _{1}\wedge \alpha _{4}+\alpha _{2}\wedge \alpha
_{6}+\alpha _{3}\wedge \alpha _{5} \\ 
\omega _{3}=\alpha _{1}\wedge \alpha _{4}+\alpha _{2}\wedge \alpha
_{6}-\alpha _{3}\wedge \alpha _{5}
\end{array}
$

\item  $\qquad 
\begin{array}[t]{l}
\left[ X_{1},X_{2}\right] =X_{6},\qquad \left[ X_{2},X_{3}\right]
=X_{5},\medskip  \\ 
\omega _{1}=\alpha _{1}\wedge \alpha _{6}+\alpha _{2}\wedge \alpha
_{5}+\alpha _{3}\wedge \alpha _{4} \\ 
\omega _{1}=-\alpha _{1}\wedge \alpha _{6}-\alpha _{2}\wedge \alpha
_{5}-\alpha _{3}\wedge \alpha _{4}
\end{array}
$

\item  $\qquad 
\begin{array}[t]{l}
\left[ X_{1},X_{2}\right] =X_{6},\medskip  \\ 
\omega =\alpha _{1}\wedge \alpha _{6}+\alpha _{2}\wedge \alpha _{5}+\alpha
_{3}\wedge \alpha _{4}
\end{array}
$

\item  $\qquad 
\begin{array}[t]{l}
\Bbb{R}^{6},\medskip  \\ 
\omega =\alpha _{1}\wedge \alpha _{6}+\alpha _{2}\wedge \alpha _{5}+\alpha
_{3}\wedge \alpha _{4}
\end{array}
$
\end{enumerate}

\end{document}